\documentclass[review]{elsarticle}

\usepackage{lineno,hyperref}
\usepackage{graphicx}
\DeclareGraphicsExtensions{.pdf,.jpeg,.png}
\usepackage{array}
\usepackage{mdwmath}
\usepackage{mdwtab}
\usepackage{fixltx2e}
\usepackage{amsmath}
\usepackage{amssymb}
\usepackage{pgfplots}
\usepackage{tikz}
\usetikzlibrary{external}
\usepackage{ifthen}
\usepackage{pdfcomment}
\usepackage[Gray,squaren,thinqspace,thinspace]{SIunits}
\usepackage{lipsum}
\usepackage{caption}
\usepackage{subcaption}
\usepackage{todonotes}
\usepackage{nicefrac}

\usepackage{color}
\definecolor{tud1d}{RGB}{36,53,114}
\definecolor{tud3d}{RGB}{0,113,94}
\definecolor{tud6d}{RGB}{174,142,0}
\definecolor{tud8d}{RGB}{169,73,19}
\definecolor{tud10d}{RGB}{115,32,84}

\usepackage[indexonlyfirst,acronym,xindy,toc,nopostdot,nomain]{glossaries}

\usepackage{amsthm}

\modulolinenumbers[5]

\journal{Computer Methods in Applied Mechanics and Engineering}

\usepackage{bm}

\newcommand{\vect}[1]{\ensuremath{\bm{\mathbf{#1}}}}										% vectors
\newcommand{\mat}[1]{\ensuremath{\mathbf{#1}}}													% matrices
\newcommand{\stiff}{\ensuremath{\mat{K}}}														    % stiffness matrix
\newcommand{\mass}{\ensuremath{\mat{M}}}															  % mass matrix
\newcommand{\coupling}{\ensuremath{\mat{G}}}													  % coupling matrix
\newcommand{\transpose}{^{\top}}																				% transpose operator
\newcommand{\parder}[2]{\ensuremath{\frac{\partial #1}{\partial #2}}}		% gradient operator

\newcommand{\grad}[1]{\ensuremath{\nabla #1}}														% gradient operator
\renewcommand{\div}[1]{\ensuremath{\nabla \cdot #1}}										% divergence operator
\newcommand{\norm}[2]{\ensuremath{\left\Vert #1 \right\Vert_{#2}}}      % norm operator
\newcommand{\Ndof}{\ensuremath{N_\text{DoF}}}    			    					   	% number of DoFs

\newcommand{\weight}{\ensuremath{\omega}}																% NURBS weights

\newcommand{\betainfsup}{\ensuremath{\beta}}														% inf-sup constant
\newcommand{\Ncoupling}{\ensuremath{N_\Gamma}}    			    						% number of harmonic coupling modes
\newcommand{\Ort}{\ensuremath{\Omega_{\mathrm{rt}}}}
\newcommand{\Ost}{\ensuremath{\Omega_{\mathrm{st}}}}
\newcommand{\Oq}{\ensuremath{\Omega_q}}
\newcommand{\oOrt}{\ensuremath{\overline{\Omega}_{\mathrm{rt}}}}
\newcommand{\oOst}{\ensuremath{\overline{\Omega}_{\mathrm{st}}}}
\newcommand{\Gd}{\ensuremath{\Gamma_{\mathrm{d}}}}
\newcommand{\Gl}{\ensuremath{\Gamma_{\mathrm{l}}}}
\newcommand{\Gr}{\ensuremath{\Gamma_{\mathrm{r}}}}
\newcommand{\Gag}{\ensuremath{\Gamma_{\text{ag}}}}
\newcommand{\nst}{\ensuremath{\vec{n}_{\mathrm{st}}}}
\newcommand{\nrt}{\ensuremath{\vec{n}_{\mathrm{rt}}}}
\newcommand{\Bst}{\ensuremath{\mat{G}_{\mathrm{st}}}}
\newcommand{\Brt}{\ensuremath{\mat{G}_{\mathrm{rt}}}}
\newcommand{\bst}{\ensuremath{\mat{g}_{\mathrm{st}}}}
\newcommand{\brt}{\ensuremath{\mat{g}_{\mathrm{rt}}}}

\newcommand{\ust}{\ensuremath{\mat{u}_{\mathrm{st}}}}
\newcommand{\urt}{\ensuremath{\mat{u}_{\mathrm{rt}}}}
\newcommand{\Kst}{\ensuremath{\mat{K}_{\mathrm{st}}}}
\newcommand{\Krt}{\ensuremath{\mat{K}_{\mathrm{rt}}}}
\newcommand{\jst}{\ensuremath{\mat{j}_{\mathrm{st}}}}
\newcommand{\jrt}{\ensuremath{\mat{j}_{\mathrm{rt}}}}
\newcommand{\Azst}{\ensuremath{A_{z,\mathrm{st}}}}
\newcommand{\Azrt}{\ensuremath{A_{z,\mathrm{rt}}}}
\newcommand{\Azq}{\ensuremath{A_{z,q}}}
\newcommand{\Htst}{\ensuremath{H_{\theta,\mathrm{st}}}}
\newcommand{\Htrt}{\ensuremath{H_{\theta,\mathrm{rt}}}}
\newacronym{pde}{PDE}{Partial Differential Equation}
\newacronym{fe}{FE}{Finite Element}
\newacronym{fem}{FEM}{Finite Element Method}
\newacronym{iga}{IGA}{Isogeometric Analysis}
\newacronym{cad}{CAD}{Computer Aided Design}
\newacronym{nurbs}{NURBS}{Non-Uniform Rational B-splines}
\newacronym{pmsm}{PMSM}{permanent magnet synchronous machine}
\newacronym{pm}{PM}{permanent magnet}
\newacronym{thd}{THD}{total harmonic distortion}
\newacronym{emf}{EMF}{electromotive force}
\newacronym{dn}{DN}{Dirichlet-to-Neumann}
\begin{document}

\begin{frontmatter}

\title{Isogeometric Analysis and Harmonic Stator-Rotor Coupling for Simulating Electric Machines}

\author[mymainaddress,mysecondaryaddress]{Zeger Bontinck \corref{mycorrespondingauthor}}
\cortext[mycorrespondingauthor]{Zeger Bontinck}
\ead{bontinck@gsc.tu-darmstadt.de}
\author[mymainaddress,mysecondaryaddress,mythirdaddress]{Jacopo Corno}
\author[mymainaddress,mysecondaryaddress]{Sebastian Sch\"ops}
\author[mysecondaryaddress]{Herbert De Gersem}

\address[mymainaddress]{Graduate School of Computational Engineering, Technische Universit\"at Darmstadt, Dolivostra\ss e 15, 64293 Darmstadt, Germany}
\address[mysecondaryaddress]{Institut f\"ur Theorie Elektromagnetischer Felder, Technische Universit\"at Darmstadt, Schlo\ss gartenstra\ss e 8, 64289 Darmstadt, Germany}
\address[mythirdaddress]{MOX Modeling and Scientific Computing, Politecnico di Milano, via Bonardi 9, 20133 Milano, Italy}

\begin{abstract}
This work proposes Isogeometric Analysis as an alternative to classical finite elements for simulating electric machines. Through the spline-based Isogeometric discretization it is possible to parametrize the circular arcs exactly, thereby avoiding any geometrical error in the representation of the air gap where a high accuracy is mandatory. To increase the generality of the method, and to allow rotation, the rotor and the stator computational domains are constructed independently as multipatch entities. The two subdomains are then coupled using harmonic basis functions at the interface which gives rise to a saddle-point problem. The properties of Isogeometric Analysis combined with harmonic stator-rotor coupling are presented. The results and performance of the new approach are compared to the ones for a classical finite element method using a permanent magnet synchronous machine as an example.
\end{abstract}


\end{frontmatter}

%\linenumbers

\section{Introduction}

\gls*{iga} was first introduced in~\cite{Hughes_2005aa,Cottrell_2009aa} and can be understood as a \gls*{fem} using a discrete function space that generalizes the classical polynomial one. \gls*{iga} has already been applied in different fields such as, e.g., mechanical engineering~\cite{Temizer_2011aa} and fluid dynamics~\cite{Gomez_2010aa}. A more elaborated overview of relevant application fields can be found in ~\cite{Nguyen_2015aa}. In this paper, we propose the application of the concepts of \gls*{iga} to electric machine simulation. According to IGA, the basis functions commonly used in \gls*{cad} for geometry construction, i.e. B-Splines and \gls*{nurbs}, are used as the basis for the solution spaces in combination with the classical \gls*{fem} framework. \gls*{iga} uses a global mapping from a reference domain to the computational domain and does not introduce a triangulation thereof. As a consequence, it is possible to represent \gls*{cad} geometries exactly, even on the coarsest level of mesh refinement. 

The possibility to parametrize circular arcs (and other conic sections) without introducing geometrical errors is of particular interest for electric machine simulation since it guarantees an exact representation of the air gap, independently of the mesh resolution. Furthermore, thanks to the properties of Isogeometric basis functions, \gls*{iga} solutions have a higher global regularity with respect to their \gls*{fem} counterparts. The inter-element smoothness of the latter is typically restricted to $C^0$. Moreover, \gls*{iga} features a better accuracy with respect to the number of degrees of freedom compared to \gls*{fem}~\cite{Cottrell_2009aa,Hughes_2014aa,Corno_2016aa}. Both advantages are of great importance for an accurate simulation of electric machines. For example, torques and forces are often calculated by the Maxwell's stress tensor evaluated in the air gap, in which case the obtained results are very sensitive to the representation and discretization of the air gap~\cite{Howe_1992aa}. This paper also tackles the problem of stator-rotor coupling which arises by our choice of \gls{iga}.

The application of \gls*{iga} for electric machine simulation is illustrated using a 2D magnetostatic formulation including the treatment of an angular displacement between stator and rotor. A further extension of the formulation to non-linear models~\cite{Pels_2015aa}, to time-harmonic \cite{Vassent_1989aa} and transient formulations~\cite{Sadowski_1993aa} and to the 3D case is straightforward.

The structure of the paper is as follows. We first introduce the 2D model commonly used to describe electric machines. We discuss how \gls*{iga} is used to discretize the model. Section~\ref{sec:coupling} first presents a naive domain decomposition approach for stator-rotor coupling and then develops a Mortar coupling strategy based on harmonic functions that is in focus of this work. A simplified example is used to analyze the convergence and stability properties of the harmonic stator-rotor coupling. Finally, we apply the proposed method for simulating a \gls*{pmsm}. The results are compared to a lowest order \gls*{fem}.

\section{IGA Electric Machine Model}
\begin{figure}
\centering
\def\svgwidth{1.0\columnwidth}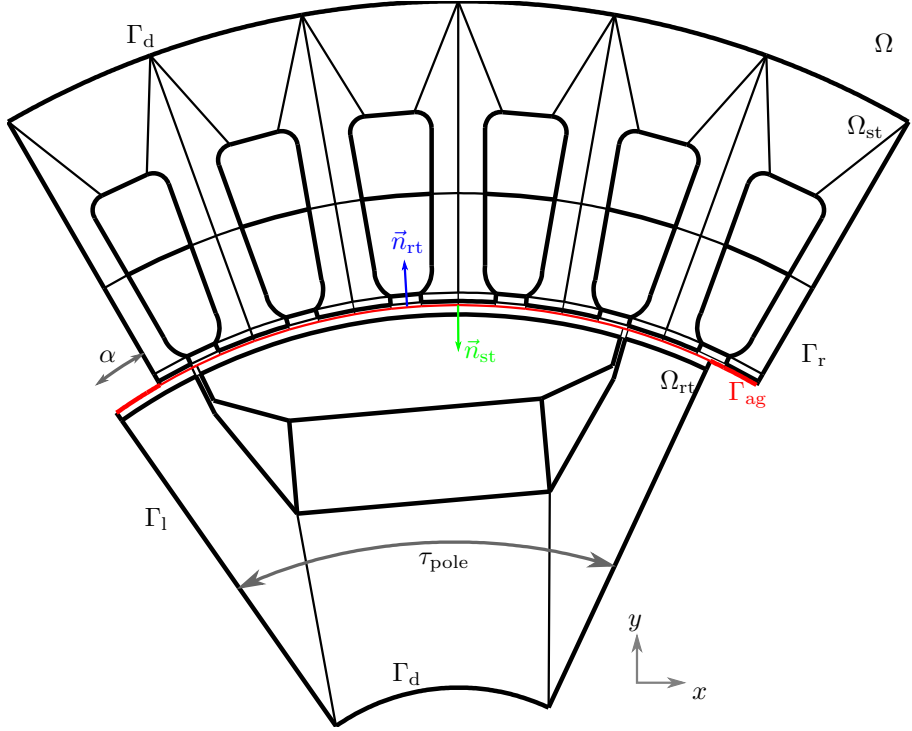
\caption{\label{fig:pole} Cross-sectional view of one pole of the machine. In red the interface between the rotor and the stator is depicted.}
\end{figure}

Electromagnetic fields are described by Maxwell's equations. For electric machines, valuable results can already be obtained using a magnetostatic formulation, i.e., a subset of Maxwell's equations where the eddy currents and displacement currents are neglected \cite{Salon_1995aa,Bianchi_2005aa}. The discretization of the resulting set of partial differential equations by FEM requires the use of N\'ed\'elec elements where the degrees of freedom are allocated on the edges of the mesh~\cite{Monk_2003aa}. Proper B-Spline approximation spaces as a counterpart to N\'ed\'elec elements in an \gls*{iga} context were introduced by Buffa et al.~\cite{Buffa_2010aa}. However, for electric machine simulation, it is often sufficient to model a 2D cross section of the geometry. Under these assumptions, Maxwell's equations reduce and combine into a  Poisson equation on the computational domain $\overline{\Omega} = \overline{\Omega}_\text{rt} \bigcup \overline{\Omega}_\text{st}$ (Fig.~\ref{fig:pole})
\begin{subequations}\label{eq:poisson}
\begin{equation}
-\nabla \cdot \left(\nu \nabla A_z\right)=\underbrace{J_{\mathrm{src}}+J_{\mathrm{pm}}}_{J_z},
\end{equation}
where $\nu=\nu(x,y)$ is the reluctivity (the inverse of the permeability), assumed to be linear and isotropic, and $A_z=A_z(x,y)$ is the $z$-component of the magnetic vector potential $\vec{A}$. The current densities exciting the coils of the machine and the magnetization current densities related to the \gls*{pm} in the rotor are depicted by $J_\mathrm{src}=J_\mathrm{src}(x,y)$ and $J_\mathrm{pm}=J_\mathrm{pm}(x,y)$, respectively. Eq.~\eqref{eq:poisson} is accompanied by Dirichlet boundary conditions at the outer stator and inner rotor boundary $\Gd$ (see Fig.~\ref{fig:pole}) and anti-periodic boundary conditions at the boundary parts $\Gl$ and $\Gr$ (see Fig.~\ref{fig:pole}), i.e.,
\begin{align}
A_z|_{\Gd} &=  0,\\
A_z|_{\Gl} &= -A_z|_{\Gr}.
\end{align}
\end{subequations}

The solution field is discretized by a linear combination of scalar basis functions $w_j(x,y)$, i.e.,
\begin{equation}\label{eq:apprx}
A_z(x,y) \approx \sum_{j=1}^{\Ndof} u_j w_j(x,y),
\end{equation}
where 
\[\mathbf{u}^\top=[u_1,\ldots,u_{\Ndof}],\]
is the vector of degrees of freedom. Applying the Galerkin approach results in the system of equations
\begin{equation}\label{eq:SoE}
\stiff\vect{u}=\vect{j},
\end{equation}
with
\begin{subequations}
\begin{equation}
k_{ij} = \int_\Omega \left(\nu\parder{w_i}{x}\parder{w_j}{x}+\nu\parder{w_i}{y}\parder{w_j}{y} \right)\;\text{d}\Omega,
\end{equation}
and writing $\mathbf{j}=\mathbf{j}_\mathrm{src}+\mathbf{j}_\mathrm{pm}$,
\begin{equation}
j_{\text{src},i} = \int_\Omega J_z w_i \;\text{d}\Omega,
\end{equation}
\begin{equation}
j_{\text{pm},i}= \int_\Omega \vec{H}_{\mathrm{pm}} \cdot \begin{bmatrix}
\parder{w_i}{y} & -\parder{w_i}{x}
\end{bmatrix}\transpose \;\text{d}\Omega.
\end{equation}
\end{subequations}
Here, $\vec{H}_{\mathrm{pm}}$ is the permanent magnet's source magnetic field strength.

There are different choices for basis functions. In this paper, two methods are considered. Firstly, there is the well established \gls*{fem} where, in the simplest case, linear hat functions are chosen \cite{Salon_1995aa}. The other approach is \gls*{iga} for which we choose \gls*{nurbs}. As the low-order \gls*{fem} can be regarded as a special case, we only discuss \gls*{iga} in the following.

\subsection{Isogeometric Analysis}

\begin{figure}
\begin{minipage}{\linewidth}
\includegraphics[width=\linewidth]{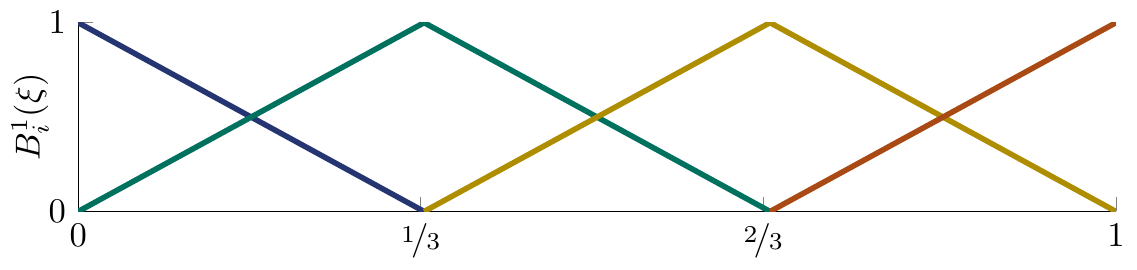}
%\begin{tikzpicture}
%	\begin{axis}[axis x line*=bottom,axis y line*=left,width=\linewidth,height=3.5cm,ylabel={$B_i^1(\xi)$},y label style={at={(axis description cs:0.1,.5)},anchor=south},xmin=0,xmax=1,ymin=0,ymax=1,xtick={0,0.33,0.66,1},
%xticklabels={0,$\nicefrac{1}{3}$, $\nicefrac{2}{3}$,1},ytick={0,1}]
%		\addplot [color=tud1d,ultra thick] table [mark=none, x=u, y=B, col sep=comma] {data/B-Spline_basis1_1.csv};
%		\addplot [color=tud3d,ultra thick] table [mark=none, x=u, y=B, col sep=comma] {data/B-Spline_basis1_2.csv};
%		\addplot [color=tud6d,ultra thick] table [mark=none, x=u, y=B, col sep=comma] {data/B-Spline_basis1_3.csv};
%		\addplot [color=tud8d,ultra thick] table [mark=none, x=u, y=B, col sep=comma] {data/B-Spline_basis1_4.csv};
%	\end{axis}
%\end{tikzpicture}
\end{minipage}
\begin{minipage}{\linewidth}
\includegraphics[width=\linewidth]{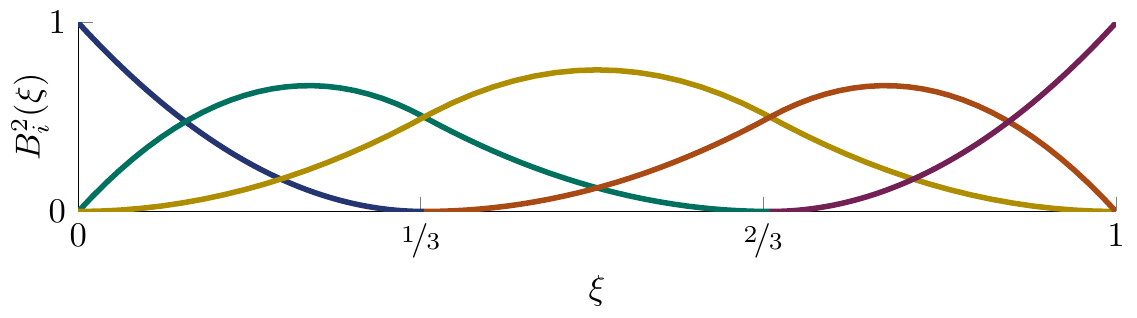}
%\begin{tikzpicture}
%	\begin{axis}[axis x line*=bottom,axis y line*=left,width=\linewidth,height=3.5cm,xlabel={$\xi$},ylabel={$B_i^2(\xi)$},y label style={at={(axis description cs:0.1,.5)},anchor=south},xmin=0,xmax=1,ymin=0,ymax=1,xtick={0,0.33,0.66,1},
%xticklabels={0,$\nicefrac{1}{3}$, $\nicefrac{2}{3}$,1},ytick={0,1}]
%		\addplot [color=tud1d,ultra thick] table [mark=none, x=u, y=B, col sep=comma] {data/B-Spline_basis2_1.csv};
%		\addplot [color=tud3d,ultra thick] table [mark=none, x=u, y=B, col sep=comma] {data/B-Spline_basis2_2.csv};
%		\addplot [color=tud6d,ultra thick] table [mark=none, x=u, y=B, col sep=comma] {data/B-Spline_basis2_3.csv};
%		\addplot [color=tud8d,ultra thick] table [mark=none, x=u, y=B, col sep=comma] {data/B-Spline_basis2_4.csv};
%		\addplot [color=tud10d,ultra thick] table [mark=none, x=u, y=B, col sep=comma] {data/B-Spline_basis2_5.csv};
%	\end{axis}
%\end{tikzpicture}
\end{minipage}
\caption{B-Spline basis functions of degree 1 and 2 on open, uniform knot vectors ($\Xi = \left[0,0,\nicefrac{1}{3},\nicefrac{2}{3},1,1\right]$ on top and $\Xi = \left[0,0,0,\nicefrac{1}{3},\nicefrac{2}{3},1,1,1\right]$ at the bottom).}\label{fig:BSp-basis-various-deg}
\end{figure}

First, we define the 1D \gls*{iga} basis functions. We choose a degree $p$ and a vector
\begin{equation}
\Xi=\begin{bmatrix}
\xi_{1} & \dots & \xi_{n+p+1}
\end{bmatrix},
\end{equation}
with $\xi_i \in \hat{\Omega} = [0,1]$, that subdivides the unit interval into elements. Here, $n$ is the dimension of the B-Spline basis which is given by the Cox-de Boor's recursion formula (see~\cite{Piegl_1997aa}). Let $\lbrace B_i^p \rbrace$ be the set of B-Spline basis functions (Fig.~\ref{fig:BSp-basis-various-deg}). We can construct the \gls*{nurbs} functions of degree $p$ as
\begin{equation}\label{eq:nurbs_basis}
N_i^p = \frac{\weight_i B_i^p}{\sum_j \weight_j B_j^p},
\end{equation}
where $\weight_i$ is a weighting factor associated to the $i$-th basis function. A general \gls*{nurbs} three-dimensional curve is obtained through the mapping
\begin{equation}\label{eq:mapping}
\mathbf{F} = \sum_{i=1}^{n} \mathbf{P}_i N_i^p,
\end{equation}
with $\mathbf{P}_i$ a set of control points in $\mathbb{R}^3$ and $\mathbf{F}:\hat{\Omega}\rightarrow\Omega\in\mathbb{R}^3$. Surfaces are built using tensor products starting from the reference square $[0,1]^2$~\cite{Piegl_1997aa}. The knot subdivision in the reference domain $\hat{\Omega}$ is transformed by the \gls*{nurbs} mapping into a physical mesh for the computational domain $\Omega$. 

Definitions~\eqref{eq:nurbs_basis} and~\eqref{eq:mapping} allow for an exact parametrization of conic sections such as circles and ellipses, which is of direct interest for the construction of the electric machine geometry, in particular of the air gap.

\gls*{iga} utilizes the same Galerkin framework as \gls*{fem}, but approximates the solution of a \gls*{pde} by a series expansion of \gls*{nurbs}~\eqref{eq:nurbs_basis}. With respect to classical \gls*{fem} spaces, the \gls*{iga} spaces bring up several advantages, e.g., the inter-element smoothness of the basis functions which allows for a higher regularity of the solution across the elements. This property also leads to a significant reduction of the number of degrees of freedom required to achieve a certain accuracy \cite{Cottrell_2009aa,Hughes_2014aa,Corno_2016aa}. The \gls*{iga} system matrices are typically smaller than their \gls*{fem} counterparts. This comes, however, at the expense of a larger bandwidth.

\section{Stator-Rotor Coupling}\label{sec:coupling} 

Due to the topology and due to the presence of different materials, it is advantageous to construct the rotor and stator model parts of the \gls*{pmsm} independently. Each part is treated by a multipatch approach, i.e., by splitting the domain into several patches, each of them being the transformation of the unit square through a \gls*{nurbs} mapping. The degrees of freedom at the interfaces between the patches are glued together through static condensation in order to obtain global $C^0$ continuity~\cite{Vazquez_2016aa}. It is worth mentioning that it would be possible to represent the entire computational domain as a single multipatch entity. Then, however, the total number of patches would increase significantly and the introduction of a relative angular displacement between stator and rotor would require a complicated mesh adaptation procedure or a complete reparametrization. 

\subsection{Domain decomposition}

In this work, we propose to subdivide the full computational domain along a circular arc $\Gag=\oOrt \bigcap \oOst$ in the air gap, separating the two subdomains $\Ort$ and $\Ost$ (Fig.~\ref{fig:pole}). The rotor and stator patches are not required to be conforming on $\Gag$. The domain-decomposition approach reads
\begin{equation}\label{eq:DD}
\left\{\begin{array}{rll}
-\nabla \cdot \left(\nu \nabla \Azq \right)& = J_{z} &\quad\text{in } \Oq, \\
\Azq|_{\Gd} &=  0, &\\
\Azq|_{\Gl} &= -\Azq|_{\Gr}, &\\
\Azrt|_{\Gag} &= \Azst|_{\Gag}, &\\
\nu \nabla\Azrt|_{\Gag}\cdot\nst &=\nu \nabla\Azst|_{\Gag}\cdot\nst, &
\end{array}\right.
\end{equation}
where $q\in\{\text{rt},\text{st}\}$ and $\nst$ is a unit vector perpendicular to the air gap interface directed from stator to rotor. The last two equations express the continuity of the magnetic vector potential and the continuity of the azimuthal component $H_\theta=\nu \nabla\Azrt|_{\Gag}\cdot\nst$ of the magnetic field strength. The two coupling approaches described below are distinct in the way they treat these interface conditions and set up the stator-rotor coupling.

\subsection{Iterative substructuring}
An iterative substructuring method invoking a Dirichlet-to-Neumann map for the first domain (here the rotor) followed by a Neumann-to-Dirichlet map for the second domain (here the stator) works as follows (see e.g.~\cite{Quarteroni_1999aa}). Let $\lambda^0$ be the initial solution for the magnetic vector potential at the air-gap interface $\Gag$ and let $k$ count the iteration steps. The domain-decomposion scheme \eqref{eq:DD} is carried out iteratively. First, a Poisson problem is solved for the rotor taking $\gamma^k$ as Dirichlet boundary condition at $\Gag$, i.e.,
\begin{equation}\label{eq:DtN-Rotor}
\left\{\begin{array}{rll}
-\nabla \cdot \left(\nu \nabla \Azrt^{k+1} \right)& = J_{z} &\text{in }\Ort, \\
\Azrt^{k+1}|_{\Gd} &=  0, &\\
\Azrt^{k+1}|_{\Gl} &= -\Azrt^{k+1}|_{\Gr}, &\\
\Azrt^{k+1}|_{\Gag} &= \gamma^k, &
\end{array}\right.
\end{equation}
from which the Neumann data $\nu\nabla\Azrt^{k+1}|_{\Gag}\cdot\nst$ is derived. Then, a Poisson problem is solved for the stator enforcing this Neumann data at $\Gag$, i.e.,
\begin{equation}\label{eq:NtD-Stator}
\left\{\begin{array}{rll}
-\nabla \cdot \left(\nu \nabla \Azst^{k+1} \right)& = J_{z} &\text{in } \Ost, \\
\Azst^{k+1}|_{\Gd} &=  0, &\\
\Azst^{k+1}|_{\Gl} &= -\Azst^{k+1}|_{\Gr}, &\\
\nu \nabla\Azst^{k+1}|_{\Gag}\cdot\nst &=
\nu \nabla\Azrt^{k+1}|_{\Gag}\cdot\nst, &
\end{array}\right.
\end{equation}
from which updated Dirichlet data
\begin{equation}
\lambda^{k+1} = \alpha \Azst^{k+1} + (1-\alpha)\gamma^k,
\end{equation}
with $\alpha\in[0,1]$ a relaxation parameter, is obtained. A relaxation factor $\alpha<1$ is required to guarantee the convergence of the iterative substructuring approach~\cite{Quarteroni_1999aa}. As a stopping criterion for the method, the $L^2$ errors for the stator and rotor models between two successive iterations should be below a user-definied tolerance, both in the rotor and the stator, i.e.,
\begin{align*}
\varepsilon_\text{rt} = \norm{\Azrt^{k+1} - \Azrt^k}{L^2(\Ort)} / \norm{\Azrt^{k+1}}{L^2(\Ort)} &< \texttt{tol}, \quad\text{and} \\
\varepsilon_\text{st} = \norm{\Azst^{k+1} - \Azst^k}{L^2(\Ost)} / \norm{\Azst^{k+1}}{L^2(\Ost)} &< \texttt{tol}.
\end{align*}

The discretization of problems~\eqref{eq:DtN-Rotor}-\eqref{eq:NtD-Stator} is carried out in the \gls*{iga} framework presented above. More details on this particular coupling approach can be found in~\cite{Bhat_2017aa}. In the following, a more sophisticated method is proposed.

\subsection{Harmonic stator-rotor coupling}
\label{sec:coupling_theory}
Let the polar coordinate system $(r,\theta)$ be connected to the stator and the polar coordinate system $(r,\theta')$ be connected to the rotor. We denote by $\alpha$ the angular displacement between the two domains, i.e. $\theta' = \theta - \alpha$. The interface conditions at $\Gag$ read
\begin{equation}
	\left\{\begin{array}{rl}
		\Azst|_{\Gag}(\theta) &= \Azrt|_{\Gag}(\theta-\alpha), \\
		\Htst(\theta) &= \Htrt(\theta-\alpha),
	\end{array}\right.
\end{equation}
where $\Htst=\nu \nabla\Azst|_{\Gag}\cdot\nst$ and $\Htrt=\nu \nabla\Azrt|_{\Gag}\cdot\nst$.

The idea of harmonic stator-rotor coupling \cite{De-Gersem_2004ad} is to express $\Htst(\theta)$ and $\Htrt(\theta')$ in terms of a particular choice of basis functions. This approach can be interpreted in the context of Mortaring methods, with a particular choice of the space of Lagrange multipliers. A superposition of harmonic functions yields
\begin{align}
	\Htst(\theta)  &=\sum_{\ell\in L} \lambda_{\mathrm{st},\ell} e^{-\imath\ell\theta} ,\\
	\Htrt(\theta') &=\sum_{\ell\in L} \lambda_{\mathrm{rt},\ell} e^{-\imath\ell\theta'} ,
\end{align}
where $\lambda_{\mathrm{st},\ell}$  and $\lambda_{\mathrm{rt},\ell}$ are the Fourier coefficients acting as degrees of freedom at the interface and $L$ is a given set of $\Ncoupling$ harmonics. The set $L$ only contains harmonic orders $\ell$ for which the corresponding harmonic functions $e^{-\imath\ell\theta}$ and $e^{-\imath\ell\theta'}$ fulfill the same anti-periodic boundary conditions as applied to $\Gl$ and $\Gr$, e.g., $e^{-\imath\ell\theta}=e^{-\imath\ell(\theta+\tau_\text{pole})}$, where $\tau_\text{pole}$ is the angular extend of one pole (Fig.~\ref{fig:pole}). The choice of harmonic trial functions enables us to construct a conforming discretization for $\Htst$ and $\Htrt$ at $\Gag$ and facilitates the application of the tangential continuity of the magnetic field strength in a strong way, leading to
\begin{equation}\label{eq:lb}
	\lambda_{\mathrm{st},\ell} =r_{\ell\ell}(\alpha)\lambda_{\mathrm{rt},\ell},
\end{equation}
where the phase shifts $r_{\ell\ell}(\alpha)$ are gathered in the rotation matrix $\mathbf{R}(\alpha)$ such that \eqref{eq:lb} is in short $\bm{\lambda}_\mathrm{st}=\mathbf{R}(\alpha)\bm{\lambda}_\mathrm{rt}$.

The discretization of the Poisson equation in $\Ost$ and $\Ort$ along the IGA framework leads to
\begin{align}
\label{eq:SoE_st}
\Kst\ust+\bst &=\jst, \\
\label{eq:SoE_rt}
\Krt\urt+\brt &=\jrt,
\end{align}
where the additional term $\bst$ and $\brt$ follow from integration by parts and are given by
\begin{align}
g_{\mathrm{st},i} &=-\int_{\Gag}\Htst(\theta) w_i(\theta)\;\text{d}\theta, \\
g_{\mathrm{rt},i} &=\int_{\Gag}\Htrt(\theta') w_i(\theta')\;\text{d}\theta'.
\end{align}
The introduction of the discretization of $\Htst$ and $\Htrt$ by harmonic functions leads to
\begin{align}
	\bst &=\Bst\mathbf{R}(\alpha)\bm{\lambda}_\mathrm{rt}, \\
	\brt &=\Brt\bm{\lambda}_\mathrm{rt},
\end{align}
where $\Bst$ and $\Brt \in \mathbb{R}^{\Ncoupling \times \Ndof }$ are coupling matrices containing the integrals
\begin{align}
g_{\mathrm{st},i\ell} &=-\int_{\Gag}e^{-\imath\ell\theta} w_i(\theta)\;\text{d}\theta, \\
g_{\mathrm{rt},i\ell} &= \int_{\Gag}e^{-\imath\ell\theta'}w_i(\theta')\;\text{d}\theta',
\end{align}
combining IGA basis functions and harmonic functions (in the spirit of mortaring).

The continuity of the magnetic vector potential at the air-gap interface is imposed in a weak way, i.e., using the complex conjugate of the harmonic functions as test functions. This results in
\begin{equation}
\label{eq:coupling}
-\mathbf{R}(\alpha)\Bst^H\ust+\Brt^H\urt=0,
\end{equation}
where this expression turns out to include the hermitian transposes of the already calculated matrices.

Combining \eqref{eq:SoE_st}, \eqref{eq:SoE_rt} and \eqref{eq:coupling} leads to the saddle-point problem
\begin{equation}
\label{eq:SoE_couple}
\left[
\begin{array}{cc}
\mathbf{K} & \mat{G}^{H}(\alpha)\\
\mat{G}(\alpha) & 0 
\end{array}\right] 
\left[\begin{array}{c}
\mathbf{U}\\
\bm{\lambda}_\text{rt}
\end{array}\right]=
\left[\begin{array}{c}
\mathbf{J}\\
0
\end{array}\right],
\end{equation}
with the blocks
\begin{equation*}
\mathbf{K}= \begin{bmatrix}
\Kst & 0 \\
0 & \Krt
\end{bmatrix},\;
\mat{G}(\alpha)= \begin{bmatrix}
-\mathbf{R}(\alpha)\Bst & \Brt
\end{bmatrix},\;
\mat{U}= \begin{bmatrix}
\urt\\
\ust
\end{bmatrix},\;
\mat{J}= \begin{bmatrix}
\jrt\\
\jst
\end{bmatrix}.
\end{equation*}

The strategy proposed in this paper can be interpreted in the context of mortar methods \cite{Belgacem_1999aa}, where the space of Lagrange multipliers is chosen as the space spanned by harmonic functions. We refer the interested reader to \cite{Buffa_2001aa}, for the application of mortaring to non-conforming \gls{fem} discretization of electrical machines, and to \cite{Brivadis_2015aa} for an overview of Isogeometric mortar methods.

\subsection{Inf-sup condition}

Problem~\eqref{eq:SoE_couple} is a saddle-point problem and may give rise to instabilities when the number of harmonics in consideration $N_\Gamma$ is too big with respect to $N_\mathrm{DoF}$~\cite{Brezzi_2012aa}. The nature of these problems has been studied thoroughly in for example \cite{Brezzi_1974aa}. In order to guarantee stability the system should satisfy the \textit{inf-sup} condition (see e.g. \cite{Bathe_2001aa}). Here, the set $L$ only contains a few harmonic orders such that we are confident that the inf-sup condition is fulfilled in practice.

The stability of the saddle-point formulation \eqref{eq:SoE_couple} can be investigated numerically. To that purpose, we calculate the inf-sup constant $\betainfsup$ using the method proposed in~\cite{Chapelle_1993aa}. Given the stiffness matrix $\stiff$ and the coupling matrix $\coupling$, $\betainfsup$ can be estimated by solving the eigenvalue problem
\begin{equation}\label{eq:inf-sup_eigproblem}
\coupling^H\stiff^{-1}\coupling\vect{x} = \sigma^2 \mass\vect{x},
\end{equation}
with
\begin{equation}
%\mass_\Gamma = \left\lbrace m_{\Gamma,lk} \right\rbrace
m_{\ell k} = \int_{\Gag} e^{-jk\theta}e^{j\ell\theta}\;\text{d}\Gamma.
\end{equation}
Given the sequence of eigenvalues $\sigma_i$ found by solving~\eqref{eq:inf-sup_eigproblem}, the inf-sup constant can be obtained as
\begin{equation}
\betainfsup = \min_i \sigma_i.
\end{equation}
The stability of the saddle-point problem is guaranteed if $\betainfsup$ is bounded away from zero.

\subsection{Verification}\label{sec:coupling_verification}

\begin{figure}
\centering
\includegraphics[width=0.75\textwidth]{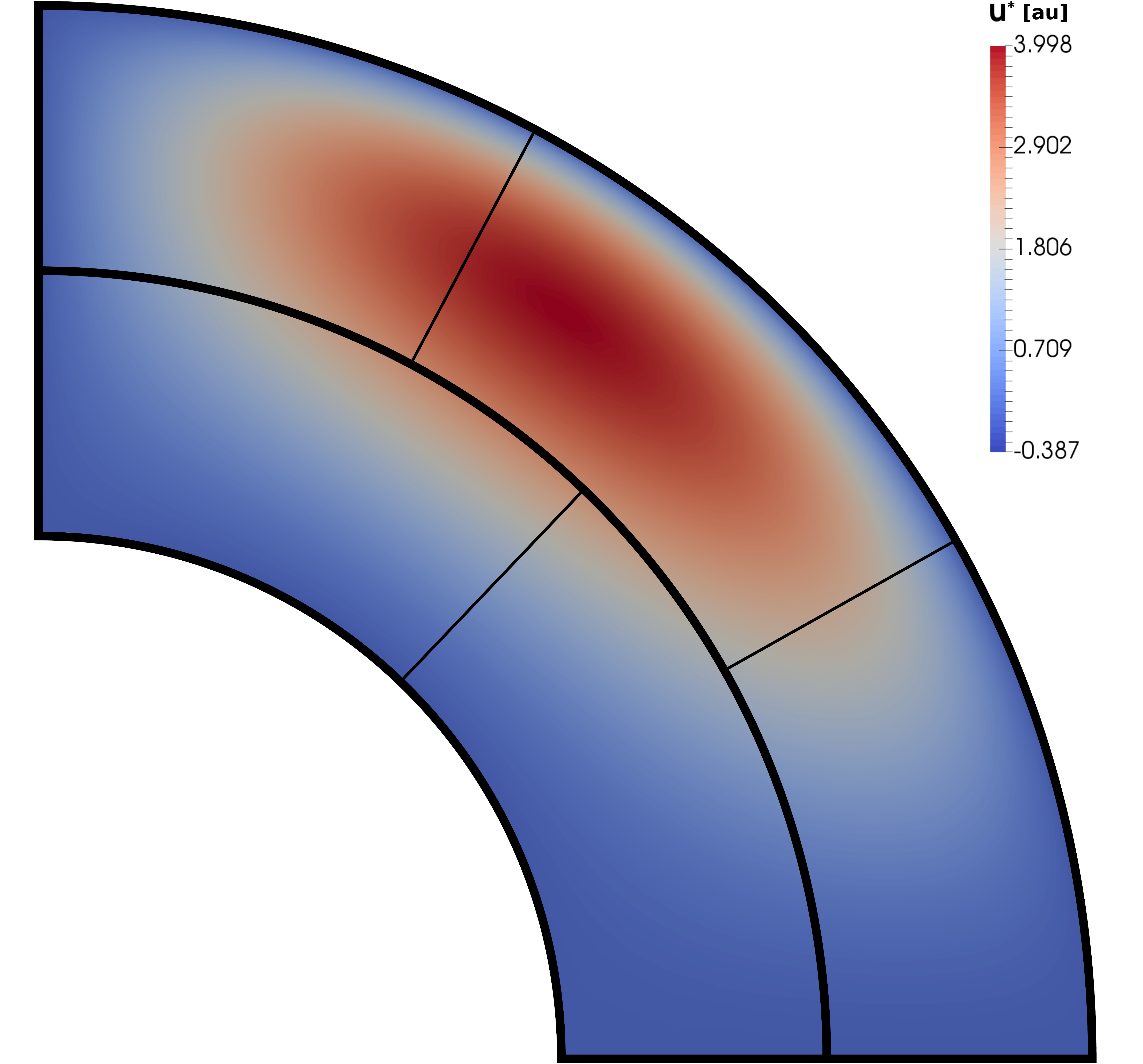}
\caption{Simplified geometry for verifying the harmonic stator-rotor coupling. The thick lines identify the two domains mimicking the rotor and the stator, while the thin lines show the non-conforming multipatch subdivisions. The color map shows the fabricated solution $u^*$.\label{fig:test_solution}}
\end{figure}

To verify the proposed harmonic stator-rotor coupling in combination with IGA, we construct a simplified example for which a closed form solution exists. The geometry chosen for this test is a quarter of a ring with an inner radius of 1 m and an outer radius of 2 m (see Fig.~\ref{fig:test_solution}). It is split in two annular domains to mimic the machine. Each domain is constructed as a multipatch geometry in such a way that the patch faces do not match at the connecting interface $\Gag$. To keep the analogy with the machine case, we will denote by $\Ort$ and $\Ost$ the inner and outer domain respectively and $\overline{\Omega}=\overline{\Omega}_\text{rt} \bigcup \overline{\Omega}_\text{st}$. As a final simplification, we consider the Poisson's equation with homogeneous Dirichlet boundary conditions
\begin{equation}
\label{eq:pois_test}
\left\lbrace
\begin{array}{rllrl}
-\div{\left(\grad{u}\right)} &= f && \text{in }\Omega,\\
u             & = 0 && \text{on }\partial\Omega,
\end{array}
\right.
\end{equation}
where 
\[f= 2x\left(22x^2y^2+21y^4-45y^2+x^4-5x^2+4\right).\]
This source term is chosen such that a closed form solution can be calculated~\cite{Vazquez_2016aa}.

\begin{figure}
\centering
\includegraphics[width=.75\linewidth]{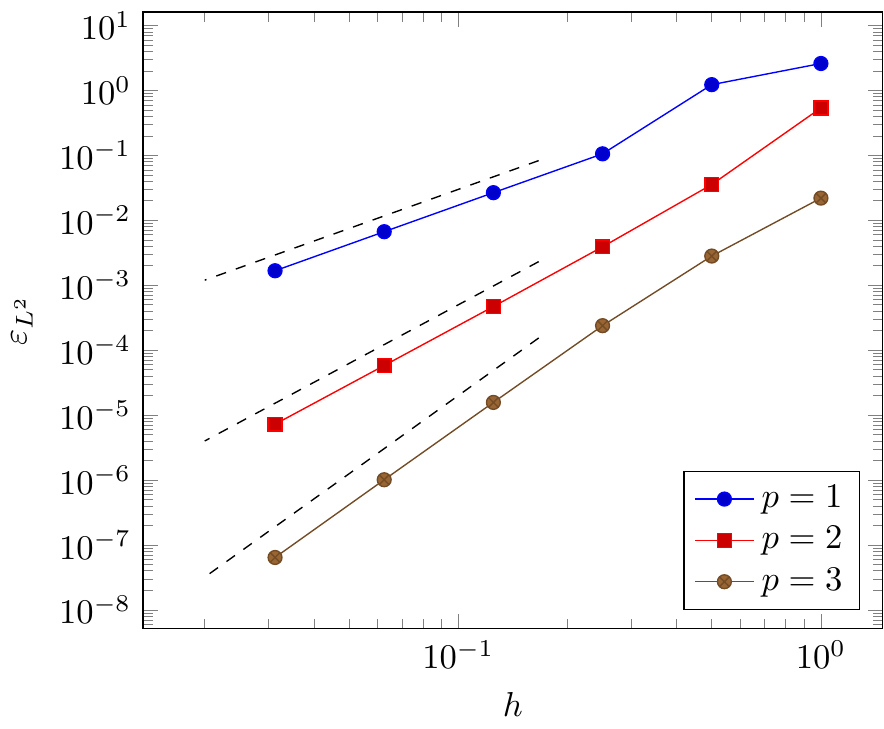}
%\begin{tikzpicture}
%\begin{loglogaxis}[width=0.75\columnwidth,xlabel={$h$},ylabel={$\varepsilon_{L^{2}}$},legend pos=south east]
%  \addplot table [x={h}, y={errl2deg1}, col sep=comma] {data/test_coupling_harmonics7.csv};
%  \addplot table [x={h}, y={errl2deg2}, col sep=comma] {data/test_coupling_harmonics7.csv};
%  \addplot table [x={h}, y={errl2deg3}, col sep=comma] {data/test_coupling_harmonics7.csv};
%  \legend{$p=1$,$p=2$,$p=3$}
%  \addplot [color=black, dashed] table [x={h}, y={order2}, col sep=comma] {data/test_coupling_harmonics7_order.csv};
%  \addplot [color=black, dashed] table [x={h}, y={order3}, col sep=comma] {data/test_coupling_harmonics7_order.csv};
%  \addplot [color=black, dashed] table [x={h}, y={order4}, col sep=comma] {data/test_coupling_harmonics7_order.csv};
%\end{loglogaxis}
%\end{tikzpicture}
\vspace{-1.5em}
\caption{Convergence of the computed solution to the exact solution for different choices of the discretization degree and for increasing mesh refinement. $\Ncoupling=7$ harmonics are used for the coupling ($\ell = -3, \dots, 3$). The dashed lines show the expected order of convergence.\label{fig:test_errl2_harmonics7}}
\end{figure}

The coupling presented in section~\ref{sec:coupling_theory} is tested for different choices of discretization degrees and number of coupling harmonics $\Ncoupling$. In Fig.~\ref{fig:test_errl2_harmonics7} the $L^2$ error with respect to the exact solution is depicted for degree $p=1,2,3$ and $\Ncoupling=7$, with $L = \lbrace -3, -2, \dots 3 \rbrace$. We recall that the $L^2$ error is defined as
\begin{equation}
\varepsilon_{L^{2}}:=\norm{u - u^*}{L^{2}} = \sqrt{\int_\Omega \left(u - u^*\right)^2 \;\text{d}\Omega}.
\end{equation}
As seen in Fig.~\ref{fig:test_errl2_harmonics7}, the coupling does not hinder the expected order of convergence, i.e. $O(h^{p+1})$.

\begin{figure}
\centering
\includegraphics[width=.75\linewidth]{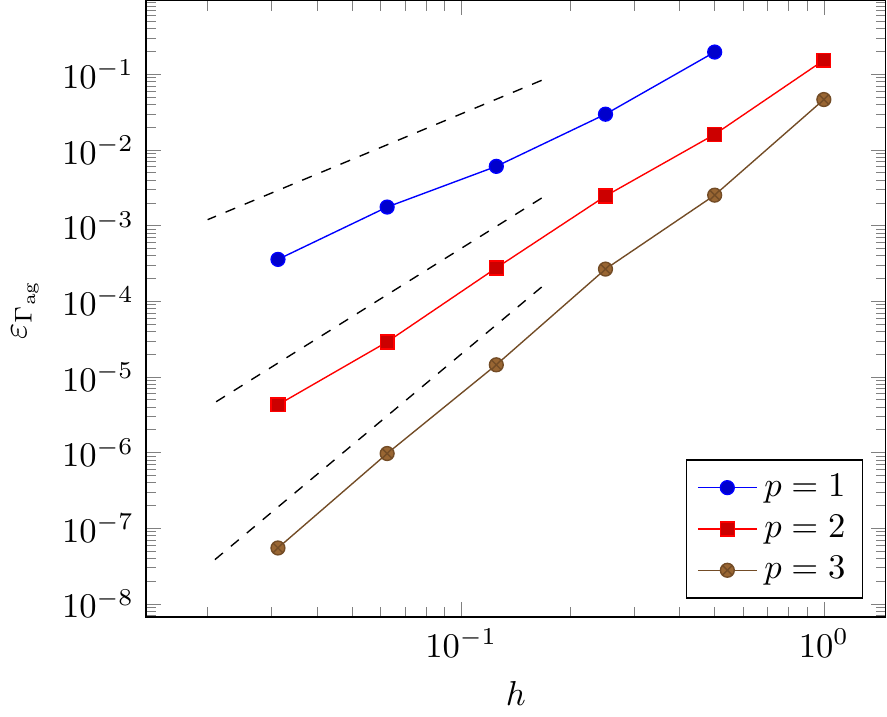}
%\begin{tikzpicture}
%\begin{loglogaxis}[width=0.75\columnwidth,xlabel={$h$},ylabel={$\varepsilon_{\Gamma_\mathrm{ag}}$},legend pos=south east]
%  \addplot table [x={h}, y={jumpl2deg1}, col sep=comma] {data/test_coupling_harmonics7.csv};
%  \addplot table [x={h}, y={jumpl2deg2}, col sep=comma] {data/test_coupling_harmonics7.csv};
%  \addplot table [x={h}, y={jumpl2deg3}, col sep=comma] {data/test_coupling_harmonics7.csv};
%	\legend{$p=1$,$p=2$,$p=3$}
%  \addplot [color=black, dashed] table [x={h}, y={order2}, col sep=comma] {data/test_coupling_harmonics7_order.csv};
%  \addplot [color=black, dashed] table [x={h}, y={order3}, col sep=comma] {data/test_coupling_harmonics7_order.csv};
%  \addplot [color=black, dashed] table [x={h}, y={order4}, col sep=comma] {data/test_coupling_harmonics7_order.csv};
%\end{loglogaxis}
%\end{tikzpicture}
\vspace{-1.5em}
\caption{Convergence of the jump of the computed solution across the interface $\Gag$ in $L^2$ norm for different choices of the discretization degree and for increasing mesh refinement. $\Ncoupling=7$ harmonics are used for the coupling ($\ell = -3, \dots, 3$). The dashed lines show the expected order of convergence.\label{fig:test_errinfty_harmonics7}}
\end{figure}

The proposed coupling imposes weak continuity of $u$ across $\Gag$. As a further verification, we compute the jump of the computed solution across $\Gag$ and evaluate its $L^2$ norm, i.e.,
\begin{equation}
\varepsilon_{\Gamma_\mathrm{ag}}:=\norm{u_{\text{rt}\vert\Gag} - u_{\text{st}\vert\Gag}}{L^{2}}.
\end{equation}
The results, depicted in Fig.~\ref{fig:test_errinfty_harmonics7}, show the convergence of the method.

\begin{figure}
\centering
\includegraphics[width=.75\linewidth]{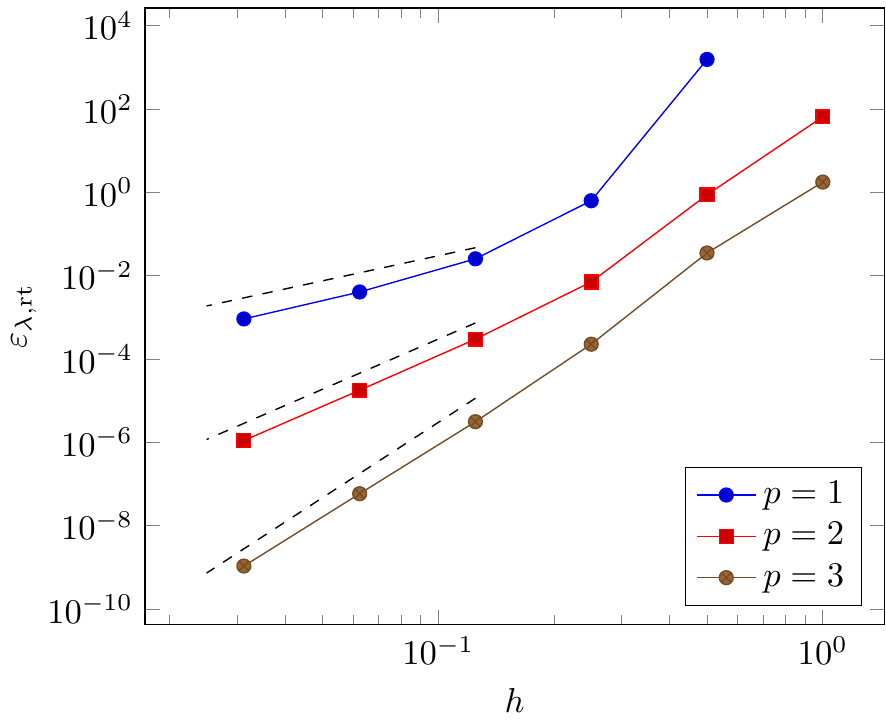}
%\begin{tikzpicture}
%\begin{loglogaxis}[width=0.75\columnwidth,xlabel={$h$},ylabel={$\varepsilon_{\lambda,\mathrm{rt}}$},legend pos=south east]
%  \addplot table [x={h}, y={errl2multipliersdeg1}, col sep=comma] {data/test_coupling_harmonics7.csv};
%  \addplot table [x={h}, y={errl2multipliersdeg2}, col sep=comma] {data/test_coupling_harmonics7.csv};
%  \addplot table [x={h}, y={errl2multipliersdeg3}, col sep=comma] {data/test_coupling_harmonics7.csv};
%	\legend{$p=1$,$p=2$,$p=3$}
%	\addplot [color=black, dashed] table [x={href}, y={order_mult_deg1}, col sep=comma] {data/test_coupling_harmonics7_mult_order.csv};
%  \addplot [color=black, dashed] table [x={href}, y={order_mult_deg2}, col sep=comma] {data/test_coupling_harmonics7_mult_order.csv};
%  \addplot [color=black, dashed] table [x={href}, y={order_mult_deg3}, col sep=comma] {data/test_coupling_harmonics7_mult_order.csv};
%\end{loglogaxis}
%\end{tikzpicture}
\vspace{-1.5em}
\caption{Convergence of the Neumann data on $\Gag$ computed from the Lagrange multipliers to the exact solution in $L^2$ norm for different choices of the discretization degree and increasing mesh refinement. $\Ncoupling=7$ harmonics are used for the coupling ($\ell = -3, \dots, 3$) The dashed lines show the order of convergence $O(h^{2p})$.\label{fig:test_errl2_multipliers_harmonics7}}
\end{figure}

Finally, we consider the convergence of the Lagrange multipliers themselves. In particular, in Fig.~\ref{fig:test_errl2_multipliers_harmonics7}, we show the convergence of the Neumann data to the exact solution, which can be evaluated as
\begin{equation}
\varepsilon_{\lambda,\mathrm{rt}}:=\norm{\grad{u^*}\cdot\vec{n}_\text{rt} - \sum_{\ell \in L} \lambda_{\text{rt},\ell} e^{-i\ell\theta}}{L^{2}}.
\end{equation}
The figure shows the good behavior of the Lagrange multipliers which converge with order $O(h^{2p})$.

\begin{figure}
\centering
\includegraphics[width=.75\linewidth]{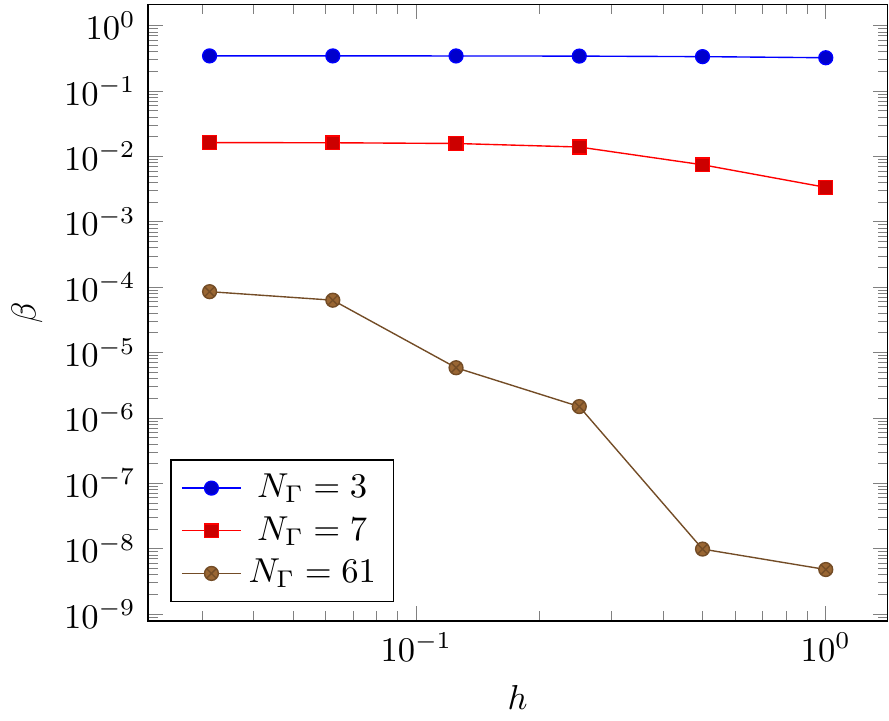}
%\begin{tikzpicture}
%\begin{loglogaxis}[width=0.75\columnwidth,xlabel={$h$},ylabel={$\betainfsup$},legend pos=south west]
%  \addplot table [x={h}, y={betadeg2}, col sep=comma] {data/test_coupling_harmonics3.csv};
%  \addplot table [x={h}, y={betadeg2}, col sep=comma] {data/test_coupling_harmonics7.csv};
%  \addplot table [x={h}, y={betadeg2}, col sep=comma] {data/test_coupling_harmonics61.csv};
%  \legend{$\Ncoupling=3$,$\Ncoupling=7$,$\Ncoupling=61$}
%\end{loglogaxis}
%\end{tikzpicture}
\vspace{-1.5em}
\caption{Numerically evaluated inf-sup constant $\betainfsup$ for a discretization of order 2 and different choices of the harmonics. When the spatial discretization is not fine enough, and when the number of harmonics increases, $\betainfsup$ goes to zero.\label{fig:test_beta}}
\end{figure}

The inf-sup constant for an Isogeometric discretization of degree $2$ and increasing number of harmonics $\Ncoupling$ is presented in Fig.~\ref{fig:test_beta}. It is apparent that when the spatial discretization is not fine enough (or, roughly said, the number of degrees of freedom is too low compared to the number of coupling harmonics), the saddle-point problem becomes unstable. For the considered applications, only a small number of harmonics with low orders is relevant. Expert knowledge can be used to choose the set $L$ appropriately.

\newpage

\section{Application: Permanent magnet synchronous machine}

IGA with harmonic stator-rotor coupling is illustrated for a $6$-pole \gls*{pmsm}. Its design is described in \cite{Henneberger_1997aa}. The longitudinal length of the machine is 10~cm. The machine is constructed of laminated steel, which is modeled as iron with a vanishing conductivity. The \gls*{pmsm} is equipped with a $3$-phase, $6$-pole distributed double-layer winding with $12$ turns per half slot.  The rotor contains $6$ buried NdFeB-magnets. The description of the geometric parameters and the material properties can be found in \ref{sec:appendix}.

\texttt{FEMM} \cite{Meeker_2009aa} and the in-house code \texttt{Niobe} are used to solve \eqref{eq:poisson} on $\Omega$ with classical \gls*{fem}. The \gls*{iga} framework is handled by the \texttt{GeoPDEs} package \cite{Vazquez_2016aa}. In post-processing, the spectrum of the \gls*{emf}, i.e., the voltage induced in the open-circuit stator windings due to rotating the rotor at nominal speed, is calculated from the solution of \eqref{eq:poisson}. This is done under no-load conditions by applying the method proposed in \cite{Rahman_1991aa}. One of the quality features of a \gls*{pmsm} is the \gls*{thd} defined by
\begin{equation}
\label{eq:THD}
\text{THD}=\frac{\sqrt{\sum_{p=2}^\infty |E_p|^2}}{|E_1|},
\end{equation}
where $p$ represents the order and $E_p$ the Fourier coefficients of the \gls*{emf}.

\subsection{Comparison Between IGA and FEM}

The rotor and stator domain are built using 12 and 78 patches respectively (see Fig.~\ref{fig:pole}). On each subdomain, an \gls*{iga} discretization is applied and the two coupling methods introduced in section~\ref{sec:coupling} are used. 

For the iterative substructuring approach, the implementation is straightforward and care must only be taken in the portion of $\Gag$ that require anti-symmetric boundary conditions (see problem~\ref{eq:poisson}). For harmonic stator-rotor coupling, as for the test case presented in section~\ref{sec:coupling_verification}, the anti-periodic boundary conditions have an impact on the composition of the set $L$ of harmonic orders. To ensure the correct behavior at the boundaries, the harmonic orders cannot be chosen arbitrarily, i.e., they need to fulfil the anti-periodic boundary condition. The chosen set is $L_{\text{ap}}= \pm 3, \pm 9, \pm 15, \dots$, i.e., because of the $6$-pole symmetry, only multiple of $3$ need to be considered and because of the mirror symmetry within a single pole, all multiples of $6$ vanish. Furthermore, since we are looking for a real valued solution we always consider double sided spectrum, i.e. if $\ell$ is chosen, then $-\ell$ is also added to the set of coupling modes.

Table~\ref{tab:timings} reports the total computational effort of matrix assembly and solving for the two straightforward implementations. In the case of the iterative substructuring approach, the two matrices to be solved are obviously smaller than in the fully coupled system \ref{eq:SoE_couple}, but the iterative procedure leads to a slow overall procedure. In particular, the main bottleneck for the iterative substructuring approach is the evaluation of the Dirichlet and Neumann data on the common boundary $\Gag$ since, in the \gls*{iga} context, the evaluation of the solution at given points in the physical domain is very expensive. This is due to the fact that this requires computing the inverse of the \gls*{nurbs} mapping in order to obtain the corresponding points in the reference domain. Given the non-linearity of the mapping, a Newton-Raphson scheme is typically used to perform this step, which considerably slows down the algorithm.

A second advantage of the proposed harmonic stator-rotor coupling is the straightforward handling of the relative rotation: for each relative angular displacement $\alpha$ only the rotation matrix has to be re-computed, which can happen at negligible cost. With iterative substructuring, instead, the right hand sides for both problems need to be re-assembled for each $\alpha$.

\begin{table}
\centering
\caption{\label{tab:timings}Comparison of the computational efficiency between \gls{iga} with iterative substructuring and \gls{iga} with harmonic stator-rotor coupling. $t$ refers to the sum of assembly and solving time.}
\begin{tabular}{l|c||r|r|r|r||r|r}
           &     & \multicolumn{4}{c||}{\textbf{DN} (\texttt{tol}$= 10^{-3}$)} & \multicolumn{2}{c}{\textbf{Harmonic}}\\
Ref. Level & $p$ & $N_\text{DoF,rt}$ & $N_\text{DoF,st}$ & $N_\text{it}$ & Time (s) & \Ndof & Time (s)\\
\hline
 1 & 1 &   32 &   193 & 9 &   67.021 &   241 &  3.335\\
 2 & 1 &  137 &   746 & 9 &  123.462 &   908 &  6.610\\
 4 & 1 &  563 &  2932 & 9 &  248.232 &  3538 & 13.675\\
 8 & 1 & 2279 & 11624 & 9 &  528.551 & 13982 & 30.437\\
16 & 1 & 9167 & 46288 & 9 & 1200.739 & 55606 & 73.414\\
\hline
 1 & 2 &   100 &   650 &  9 &   75.401 & 771 & 3.356\\
 2 & 2 &   256 &  1515 & 10 &  161.381 & 1801 & 6.813\\
 4 & 2 &   784 &  4325 & 10 &  330.292 & 5157 & 15.004\\
 8 & 2 &  2704 & 14265 & 10 &  710.511 & 15053 & 33.953\\
16 & 2 & 10000 & 51425 &  9 & 1490.251 & 61581 & 85.595\\
\end{tabular}
\end{table}

\begin{figure}
\centering
\includegraphics[width=\linewidth]{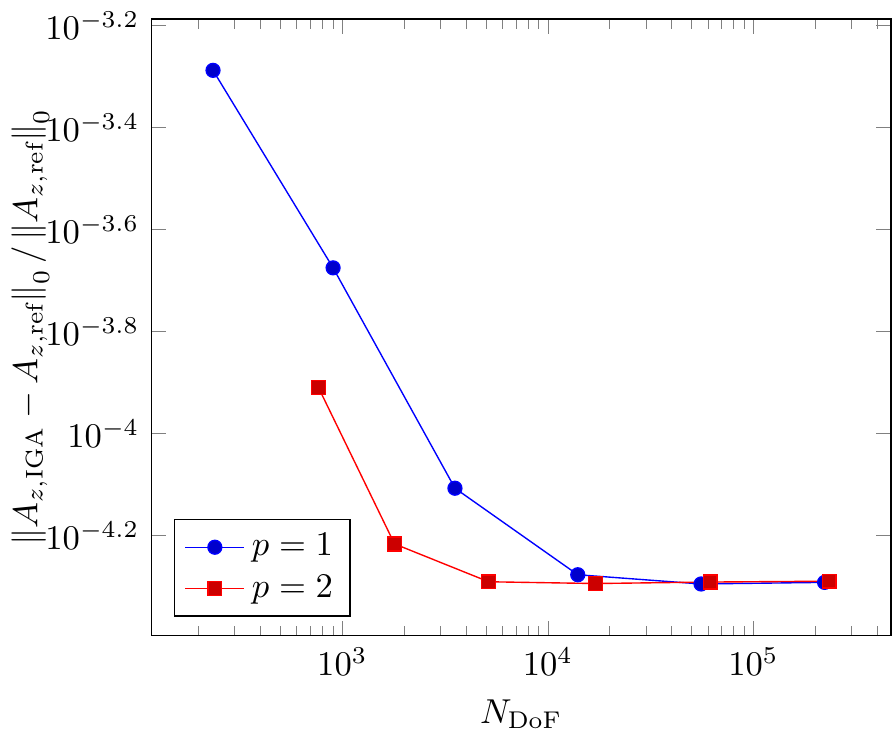}
%\begin{tikzpicture}
%\begin{loglogaxis}[width=0.75\columnwidth,xlabel={$\Ndof$},ylabel={$\norm{A_{z,\text{IGA}} - A_{z,\text{ref}}}{0} / \norm{ A_{z,\text{ref}}}{0}$},legend pos=south west]
%  \addplot table [x={ndofdeg1}, y={errl2deg1}, col sep=comma] {data/machine_pole_iga_fem_errl2.csv};
%  \addplot table [x={ndofdeg2}, y={errl2deg2}, col sep=comma] {data/machine_pole_iga_fem_errl2.csv};
%  %\addplot table [x={Ndof}, y={normsfem}, col sep=comma] {data/fem_conv.csv};
%	\legend{$p=1$,$p=2$}
%\end{loglogaxis}
%\end{tikzpicture}
\caption{Convergence of \gls*{iga} with harmonic stator-rotor coupling ($\ell = \pm3, \pm9, \pm15$) towards a first order \gls*{fem} solution on a very fine mesh ($\Ndof=225667$). \label{fig:errl2_machine}}
\end{figure}

We then consider \gls*{iga} with harmonic stator-rotor coupling and compare it to a finite element model with a conforming discretisation in the air gap. The results for discretization degree $p=1,2$ are compared to a reference solution obtained using classical first order \gls*{fem} on triangles~\cite{Salon_1995aa} in Fig.~\ref{fig:errl2_machine}. The simulation points shown in the plot correspond to the data given in Table~\ref{tab:timings}.

\begin{figure}
\centering
\includegraphics[width=\linewidth]{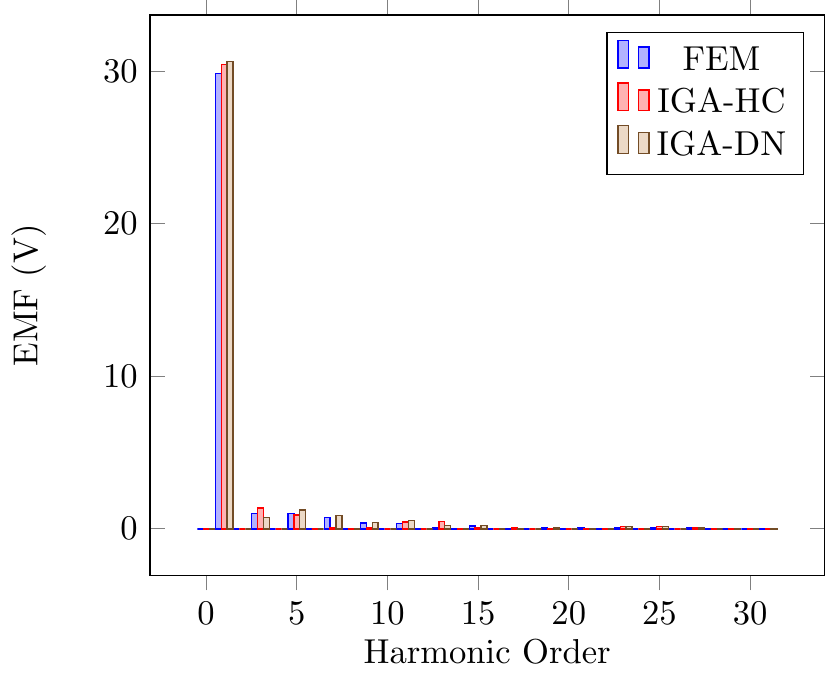}
%\begin{tikzpicture} 
%\begin{axis}[ xlabel=Harmonic Order, ylabel=\gls*{emf} (V), ybar=0pt,
%bar width=1.7pt, legend pos = north east, xtick = {0, 5, 10, 15, 20, 25, 30, 35} %, symbolic x coords={0,1,2,3,4,5,6,7,8,9,10,11,12,13,14,15,16,17,18,19,20,21,22,23,24,25,26,27,28,29,30,31}
%] 
%\addplot table [x={l_FEM}, y={E_FEM}, col sep=comma] {data/barplotdata.csv};
%\addplot table [x={l_IGA}, y={E_IGA}, col sep=comma] {data/barplotdata.csv};
%\addplot table [x={l_IGA}, y={E_IGA}, col sep=comma] {data/barplot_DN.csv};
%\legend{FEM,IGA-HC, IGA-DN} 
%\end{axis} 
%\end{tikzpicture}
\caption{Spectrum with the first 32 modes of the \gls*{emf} of the \gls*{pmsm}. \label{fig:bar}}
\label{fig:spectra}
\end{figure}

In Fig.~\ref{fig:spectra}, parts of the spectra for the EMF calculated by \gls*{fem} and by \gls*{iga} are computed. The results for the \gls*{emf} and the \gls*{thd} are listed in Table~\ref{tab:results}. The results for the \gls*{emf} have a maximal relative difference below $3\%$. The relative difference for \gls*{thd} is higher, namely $6\%$, which might imply that the iterative substructuring introduces higher harmonics.

\begin{table}
\centering
\caption{\label{tab:results}Numerical results for the \gls*{emf} and the  \gls*{thd}.}
\begin{tabular}{l|c|c|r|r}
 & $E_1$& \gls*{thd} & $\Ndof$ & Time (s)\\
\hline
FEM& 29.8 V & 5.72$\cdot$10$^{-2}$ $\%$ & 225667 & 103.45\\
IGA-HC& 30.4 V & 5.87$\cdot$10$^{-2}$ $\%$ & 5157 & 15.00\\
IGA-DN& 30.6 V & 6.06$\cdot$10$^{-2}$ $\%$ & 5109 & 330.29\\
%\hline
%Relative Difference& 1.8 $\%$ & 2.8 $\%$
\end{tabular}
\end{table}

\section{Conclusion}
In this work \gls{iga} has been applied to model a \gls{pmsm}. Since it is possible to parametrize circular arcs exactly, geometric approximations, from which the classical \glspl{fem} suffer, are avoided. A multipatch approach is used to model the rotor and the stator separately. The coupling between the two parts has been carried out by iterative substructuring or by using harmonic basis functions. To test the latter so-called harmonic stator-rotor coupling, a test case has been constructed for which the convergence of the spatial discretization has been shown. The harmonic stator-rotor coupling leads to a saddle-point problem for which our setting attains stability. As illustrated by the example, IGA with harmonic stator-rotor coupling is a new and promising alternative to standard finite element procedures for electric machine simulation. 

\section*{Acknowledgment}
This work is supported by the German BMBF in the context of the SIMUROM project (grant nr. 05M2013), by the DFG grant \textbf{SCHO1562/3-1}, by the 'Excellence Initiative' of the German Federal and State Governments and the Graduate School of CE at TU Darmstadt. The authors would also like to thank Ms. P. Bhat for the work she did in the framework of her master thesis.
\newpage
\appendix
\setcounter{table}{0}
\setcounter{figure}{0}
\section{Machine parameters}
\label{sec:appendix}
The parameters of the machine are listed in Tab.~\ref{tab:geom_param} and in Tab.~\ref{tab:mat_param}. The geometrical parameters are depicted in Fig.~\ref{fig:machine_param}. 

\begin{figure}[h]
\centering
\def\svgwidth{0.9\columnwidth}\input{geometry_IGA_machine_2.pdf_tex}
\caption{\label{fig:machine_param} Geometry of the \gls*{pmsm}}
\end{figure}
\begin{table}
\centering

\caption{\label{tab:geom_param}Parameters describing the geometry of the machine}
\begin{tabular}{l|l|l}
\multicolumn{3}{c}{Rotor}\\
\hline
Inner radius rotor & $R_\mathrm{rt,i}$ & 16 mm\\
Outer radius rotor & $R_\mathrm{rt,o}$ & 44 mm\\
Magnet width	   & $d_1$ & 19 mm\\
Magnet height      & $d_2$ & 7 mm\\
Depth of the magnet in rotor& $d_3$   & 7 mm \\
 &$\delta_1$ & 8.5\degree\\
 &$\delta_2$ & 42\degree\\
\hline
\multicolumn{3}{c}{Stator}\\
\hline
Inner radius stator& $R_\mathrm{st,i}$ & 45 mm \\
Outer radius stator& $R_\mathrm{st,o}$ & 67.5 mm\\
Number of turns	   & $N_w$			   & 12\\
 &$\delta_3$ & 7\degree\\
 &$\delta_4$ & 5.7\degree\\
 &$\delta_5$ & 4\degree\\
 &$l_1$ & 0.6 mm\\
 &$l_2$ & 5.4 mm\\
 &$l_3$ & 5 mm\\
 &$l_4$ & 8.2 mm\\
Skew angle & $\varphi$ & 0.52\degree\\
\hline
\multicolumn{3}{c}{Air gap}\\
\hline
Radius of $\Gag$   & $R_\mathrm{ag}$   & 44.7 mm \\
\end{tabular}
\end{table}
\begin{table}
\centering

\caption{\label{tab:mat_param}Parameters describing the material properties}
\begin{tabular}{l|l|l}
\multicolumn{3}{c}{Material properties}\\
\hline
Conductivity of iron				& $\sigma_\mathrm{Fe}$ & 0 S/m \\
Conductivity of copper				& $\sigma_\mathrm{Cu}$ & 4.3$\cdot 10^7$ S/m \\
Conductivity of \gls*{pm}			& $\sigma_\mathrm{PM}$ & 6667 S/m \\
Relative permeability of iron		& $\mu_\mathrm{r,Fe}$ & 500 \\
Relative permeability of copper 	& $\mu_\mathrm{r,Cu}$ & 1 \\
Relative permeability of \gls*{pm}  & $\mu_\mathrm{r,PM}$ & 1.5 \\
Remanent magnetic field of \gls*{pm}& $B_\mathrm{r}$& 0.94 T
\end{tabular}
\end{table}
\newpage

%\printbibliography%

\end{document}